# How Does One Cut a 30-60-90 Triangle?

**Ed Pegg, Brian Trease, Andrew Hudson and D. S. McNeil**

How does one cut a 30-60-90 triangle [3] into triangles similar to itself? Pieces are unlabeled, and dissections are counted up to similarity of the entire figure, including reflection.

Let $A$, $B$, $C$ be the 90°, 30°, 60° vertices. The altitude $AH$, with foot $H$ on the hypotenuse $BC$, gives the unique two-piece dissection, with areas $1/4$ and $3/4$. Splitting either piece by its own altitude gives two three-piece dissections. For the third, bisect angle $C$, meeting $AB$ at $E$, then drop the perpendicular from $E$ to $BC$. Its three pieces have equal area.

Taking the parent triangle to have area 1, labels give the least common denominator of the piece areas, the multiset letter, and the arrangement index, as in 25e10. The numerals inside the pieces are the corresponding area numerators. For a fixed piece count, multisets are ordered by increasing denominator, then by descending lexicographic order of their numerators listed largest first. The number of pieces is a prefix when necessary, such as 4-9a1 and 5-9a1.

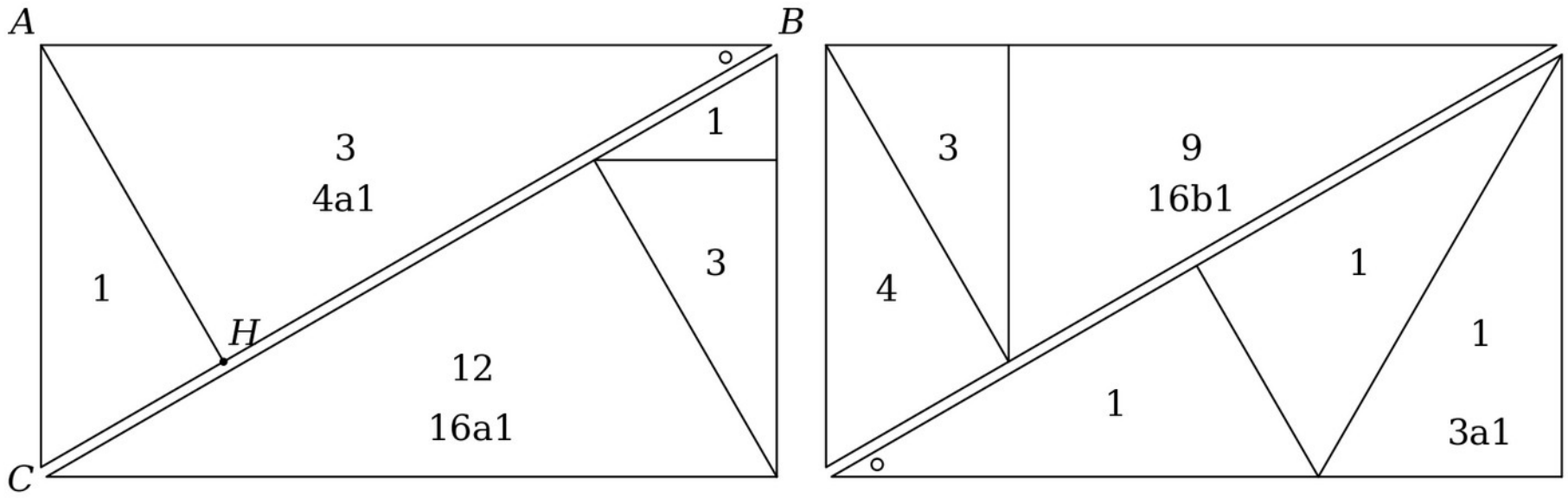

Figure 1. The two-piece dissection and the three three-piece dissections. The hypotenuse cut in 3a1 can be considered a T-junction.

## Four pieces

Figure 2 gives the eighteen four-piece arrangements. Three come from a symmedian construction. A symmedian is the reflection of a median in the corresponding angle bisector. The three symmedians meet at $K$, the symmedian point; in a right triangle, $K$ is the midpoint of $AH$ [1].

From a vertex $P$, let the symmedian meet the opposite side $QR$ at $D$. Draw $DE$ parallel to $PR$ and $DF$ parallel to $PQ$, with $E$ on $PQ$ and $F$ on $PR$, then join $EF$. The resulting four triangles are similar to the original. Choosing the 90°, 60°, and 30° vertices gives 16a2, 25a1, and 49a1, respectively. Their symmedians and the common point $K$ are marked.

A small circle inside the parent triangle’s 30° corner marks an independent dissection.

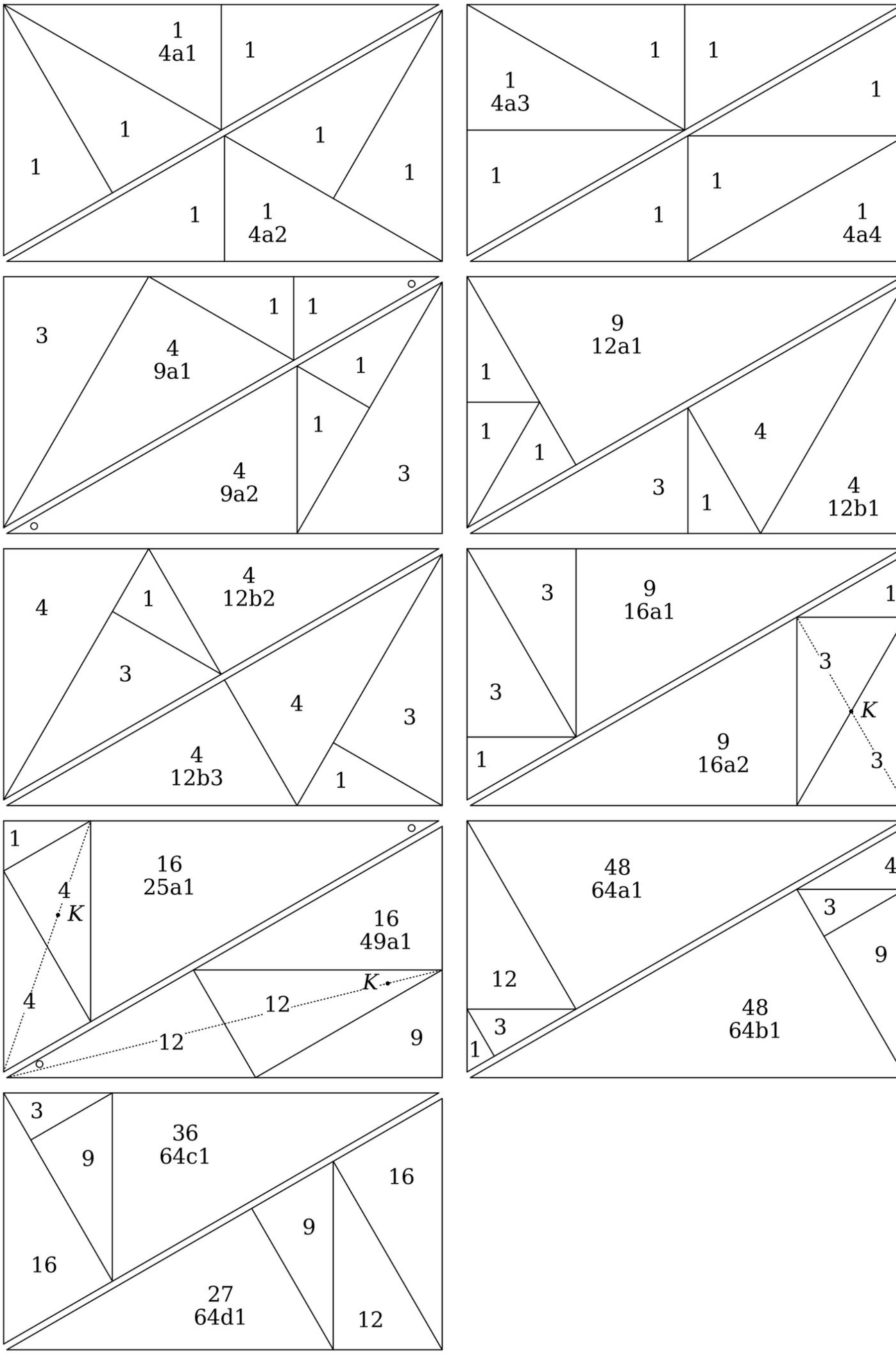

Figure 2. All eighteen four-piece arrangements. Dots mark $K$; dashed segments are construction lines.

## **Compressing by area multisets**

The eighteen arrangements use eleven distinct multisets of piece areas. Since every piece has the same angles, equal area means congruence. Figure 3 retains the last indexed arrangement in each set. Its final index therefore also counts the arrangements in that set, while every retained label refers to the identical dissection in Figure 2. Figures 4 and 6a–6c use the same convention.

For example, 4a4 is the fourth of four arrangements of congruent quarters; 12b3 is the third of three arrangements with areas $\{4,4,3,1\}/12$. The eleven final numbers sum to $4+2+1+3+2+1+1+1+1+1+1=18$.

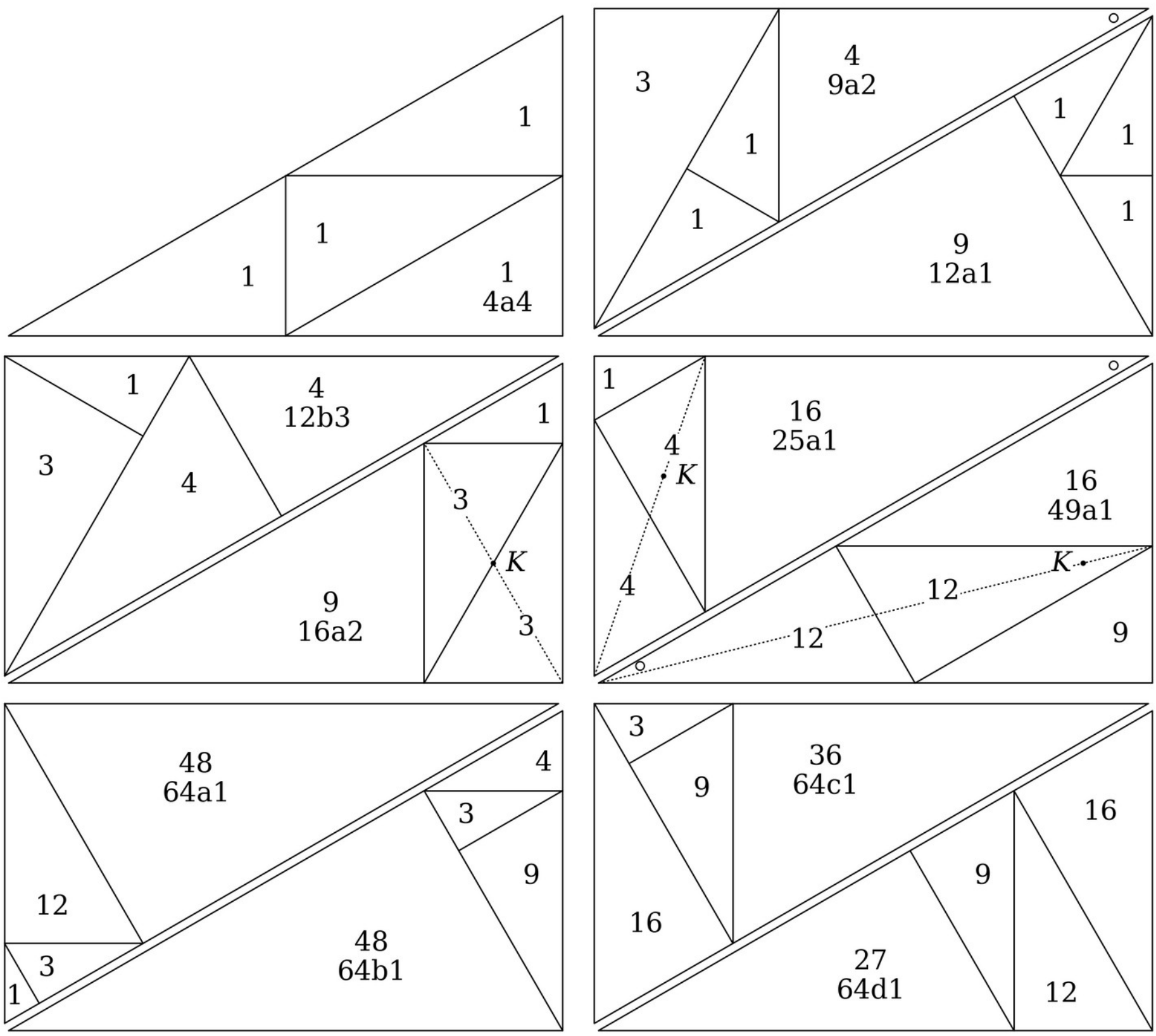

Figure 3. One representative for each area multiset, with labels and layouts matching Figure 2.

A layout is derived if it is obtained by substituting one dissection into a piece of another, both having fewer pieces than the result, or if it can be reached from such a layout by a sequence of local rotations or reflections, each preserving the union of a proper group of pieces; otherwise it is independent.

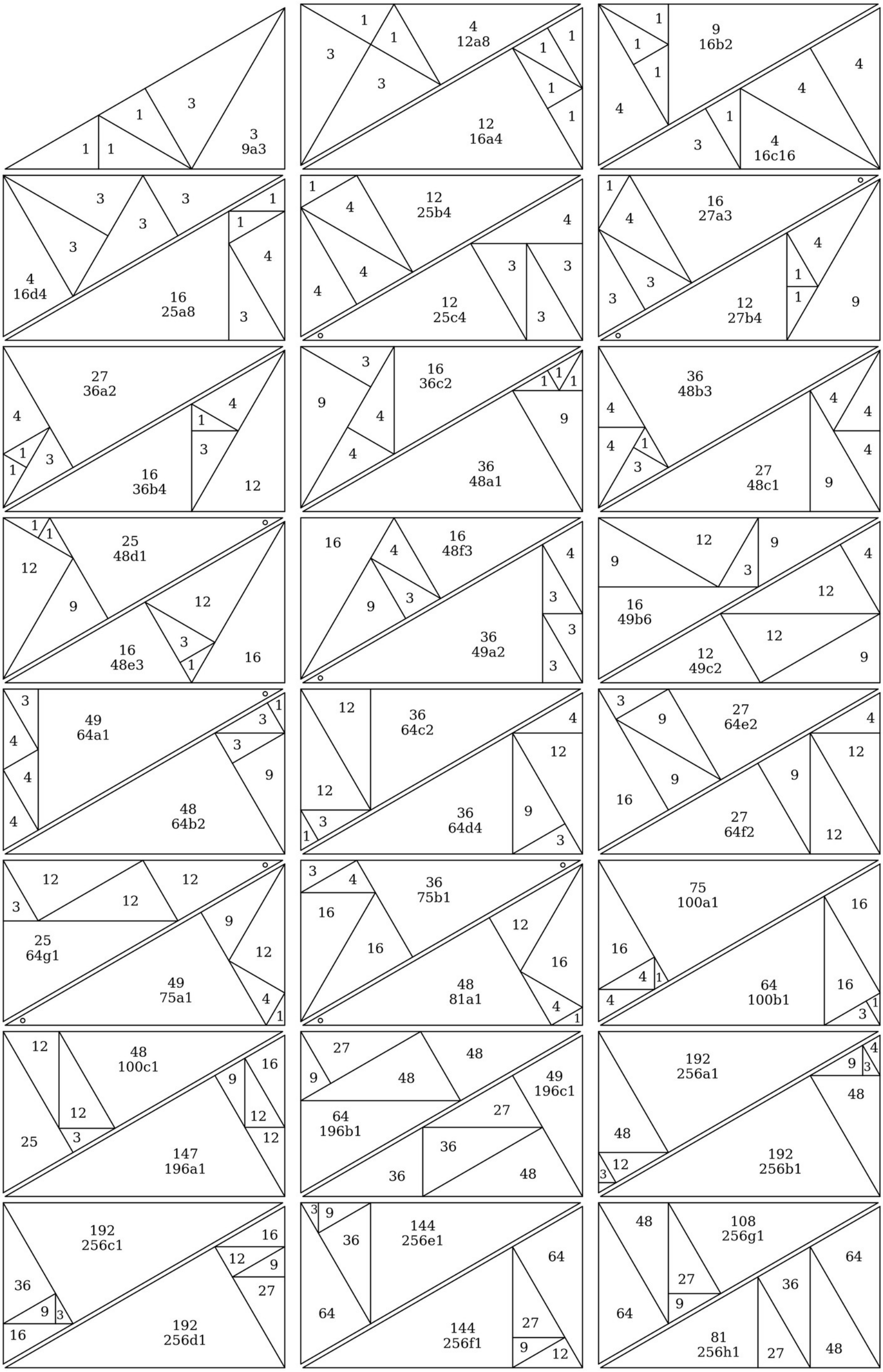

Figure 4. Forty-seven area multisets for 121 five-piece arrangements.

## Five pieces

The five-piece catalog contains 121 arrangements, using 47 distinct multisets of piece areas. Repeated subdivision accounts for 84 arrangements. A two-piece template can be inserted into a four-piece layout, a three-piece template into a three-piece layout, or a four-piece template into a two-piece layout.

The angle relations 60° = 2 × 30° and 90° = 3 × 30° permit additional corner splits. The side ratio 1:√3:2 lets one hypotenuse match two short legs of congruent pieces. Two copies joined along their long legs form an equilateral triangle, allowing local 120° rotations.

Subdivision does not produce every arrangement. A proper group of pieces can sometimes be rotated or reflected while the other pieces remain fixed. In particular, all three rotational positions of a dissection inside an equilateral region must be checked. Sequences of these operations connect eighteen further layouts to subdivision layouts.

The remaining nineteen five-piece dissections are independent, giving 84 + 18 + 19 = 121 arrangements. Figure 5 shows six examples whose area multisets each occur in only one arrangement.

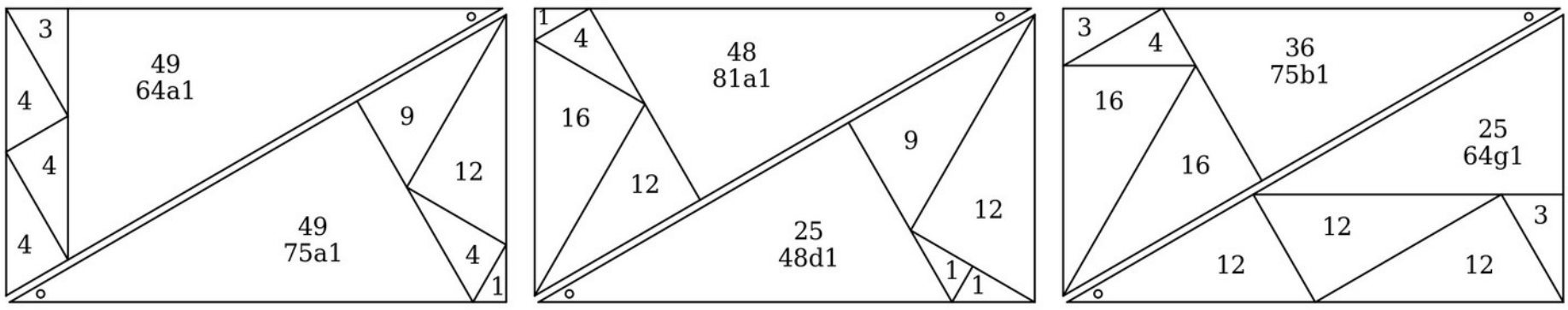

Figure 5. Six examples of independent five-piece dissections.

In 48d1 of Figure 4, the edge separating the two pieces numbered 1 ends inside a side of the piece numbered 12, forming a T-junction. Requiring every shared segment to be a complete side of both pieces would miss this layout. Allowing only T-junctions in the interior would also miss the fan in 12a8: the pieces numbered 3, 3, 1, and 1 meet at their right-angle vertices, with no triangle side passing through the junction. T-junctions also occur in larger catalogs.

## Six pieces

Figures 6a–6c show the 235 area multisets of 844 six-piece layouts. Subdivision produces 667 layouts; sequences of local rotations or reflections connect 89 more to those layouts, leaving 88 independent layouts.

Independence can depend on the arrangement. Only 9a23 is independent among the 23 arrangements of its area multiset, while the 25e multiset has two independent arrangements, including 25e10. The six-piece multiplicity record is 16b36: 36 arrangements of {4,4,3,3,1,1}/16.

The smallest wholly interior piece in the six-piece catalog has area 1/324 in 324d1, followed by 1/300 in 300f1. Both layouts are independent.

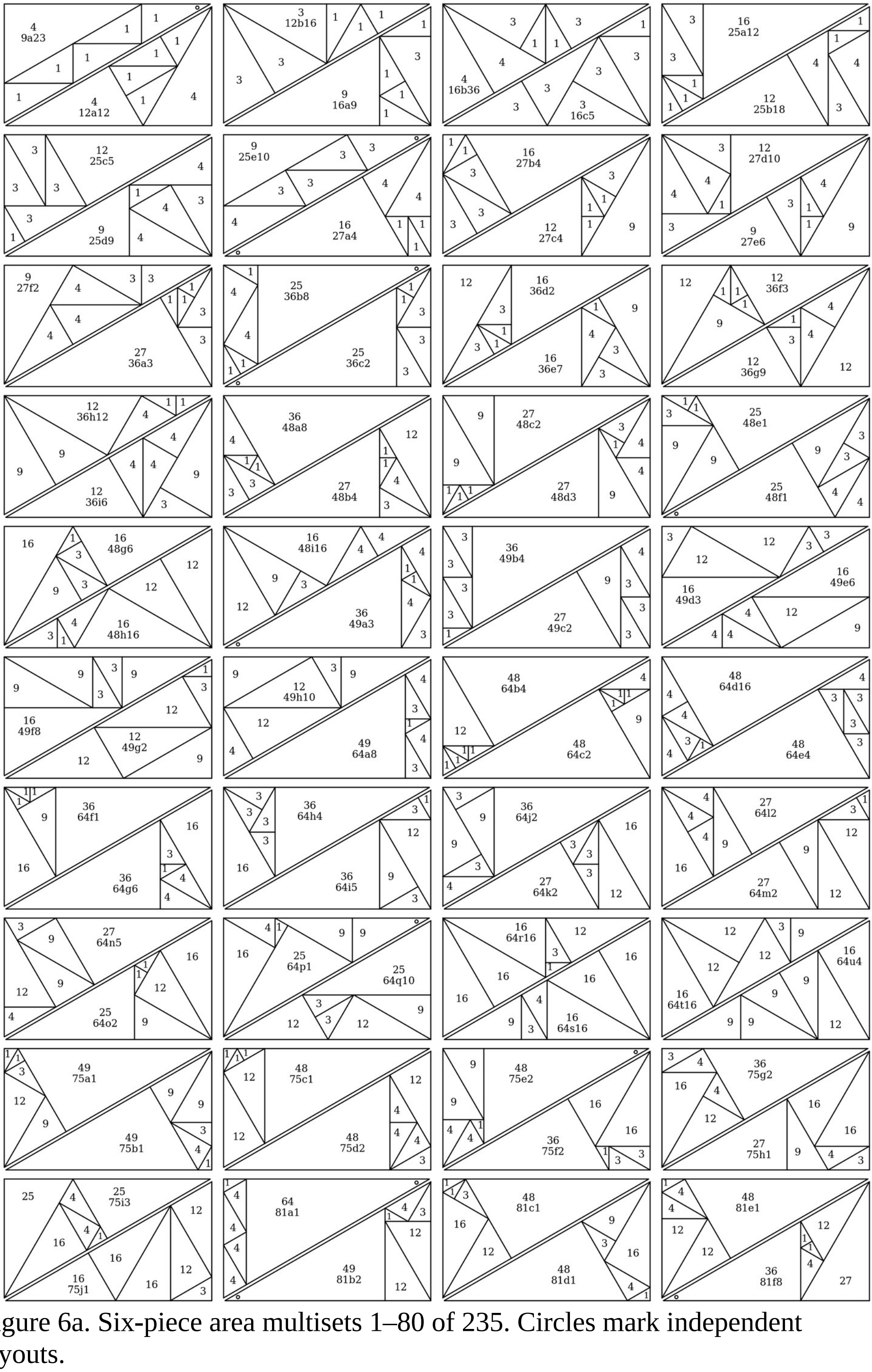

Figure 6a. Six-piece area multisets 1–80 of 235. Circles mark independent layouts.

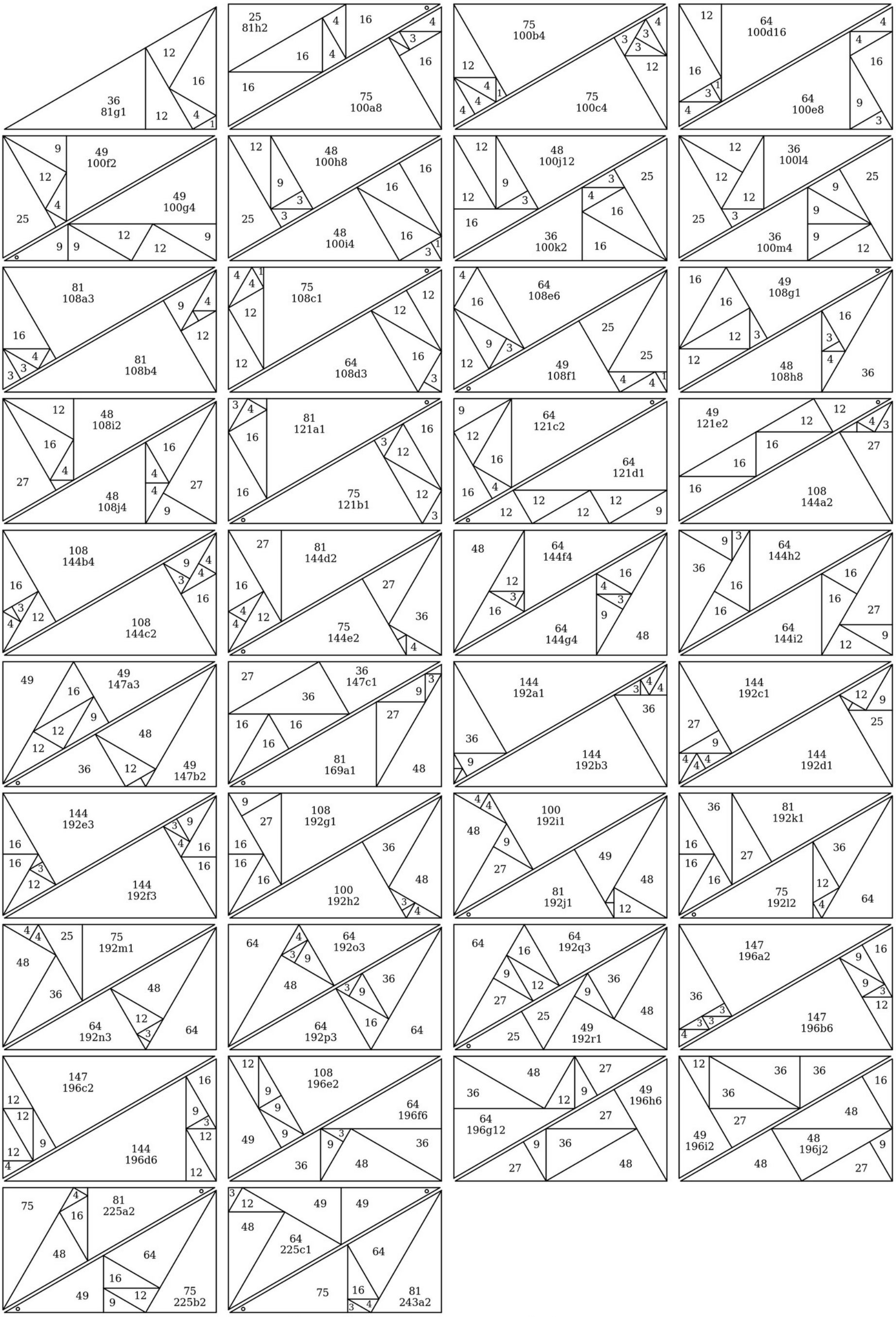

Unnumbered tiny pieces: □ 1 ■ 3

Figure 6b. Six-piece area multisets 81–155 of 235.

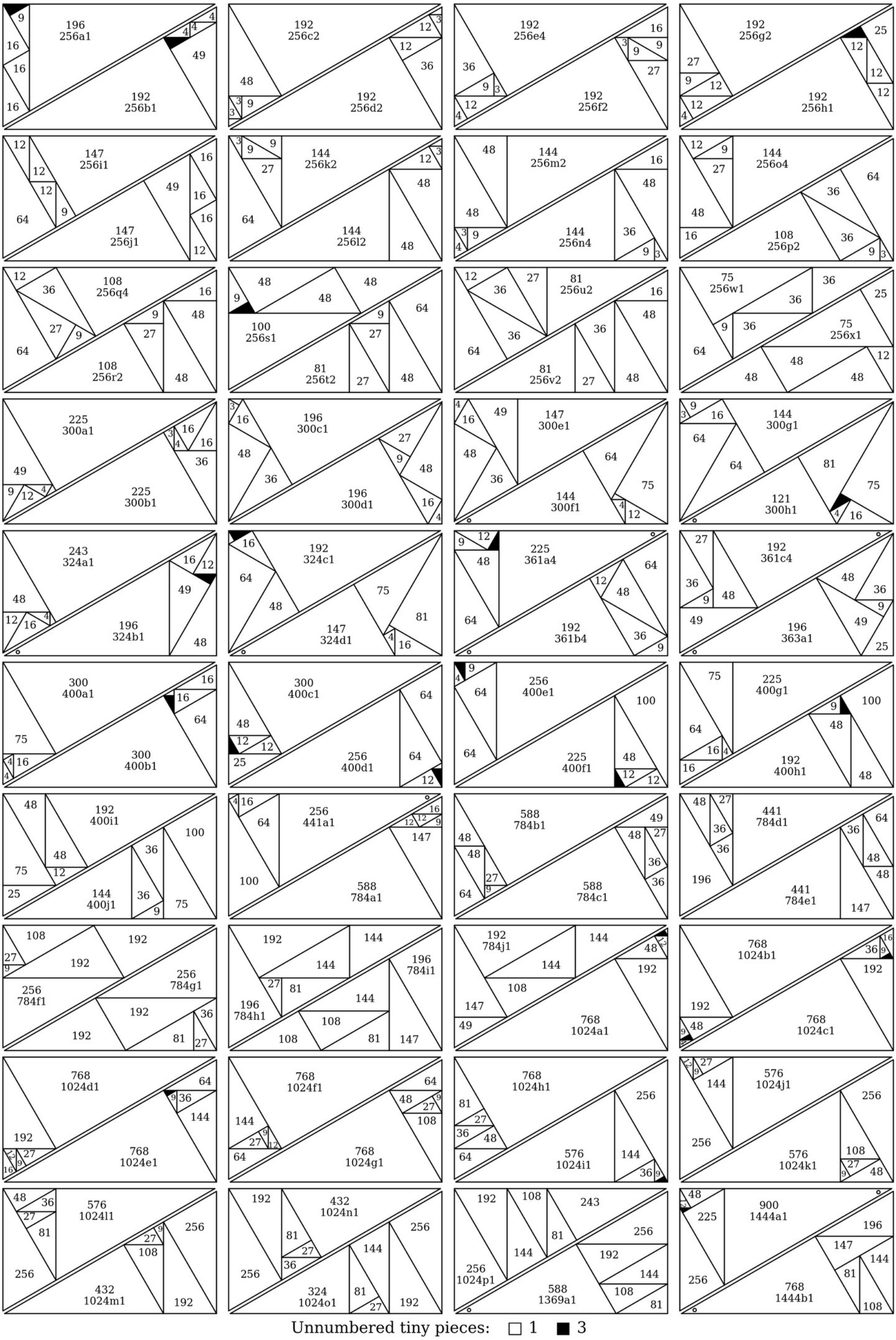

Figure 6c. Six-piece area multisets 156–235 of 235.

## Six-piece examples

Figure 7 gives six-piece examples with interior junctions: the smallest pieces in 64p1 and 192r1 touch the boundary at one vertex, while those in 192j1 and 300h1 are wholly interior.

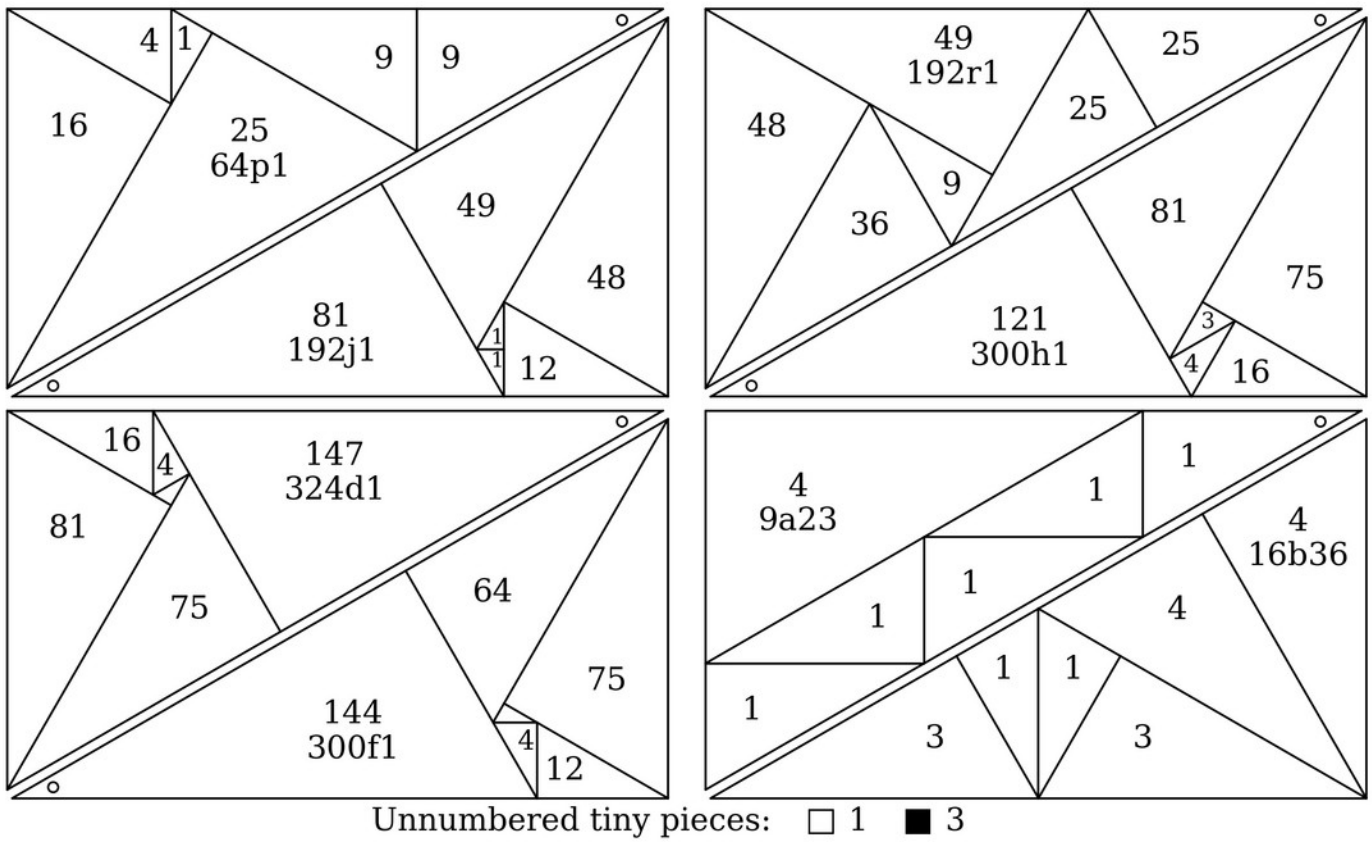

Figure 7. Selected six-piece layouts.

## Computational checks

The search and checking code were developed with assistance from ChatGPT (OpenAI). Removing the piece at the unsplit outer 30° angle leaves a triangle or convex quadrilateral. The search varies boundary cut points and interior junctions together. A separately implemented search reached nine pieces, finding 460404 arrangements at nine and no unresolved singular systems. It matched every layout in the earlier catalog through eight pieces.

Right, non-right and isosceles denote generic triangles.

| Triangle | OEIS | 1 | 2 | 3 | 4 | 5 | 6 | 7 | 8 | 9 |
|---|---|---|---|---|---|---|---|---|---|---|
| 30-60-90 | A399999 | 1 | 1 | 3 | 18 | 121 | 844 | 6397 | 53145 | 460404 |
| Right | A399770 | 1 | 1 | 2 | 10 | 60 | 363 | 2346 | 16581 | 123182 |
| 45-45-90 | A400511 | 1 | 1 | 1 | 5 | 20 | 100 | 543 | 3209 | 19986 |
| Non-right | A399765 | 1 | 0 | 0 | 4 | 0 | 36 | 64 | 468 | 1792 |
| Isosceles | A400510 | 1 | 0 | 0 | 2 | 0 | 11 | 19 | 108 | 385 |
| Equilateral | A299705 | 1 | 0 | 0 | 1 | 0 | 1 | 2 | 3 | 9 |

Sequence A399999 (formerly A400462) was a runner-up in the 2026 OEIS 400K celebration [4].